\documentclass[11pt]{amsart}
\usepackage{amsmath,amssymb, xypic,verbatim,amscd}

\swapnumbers
\theoremstyle{plain}
\newtheorem{theorem}{Theorem}[section]

\newtheorem{lemma}[theorem]{Lemma}

\newtheorem{proposition}[theorem]{Proposition}

\theoremstyle{definition}
\newtheorem{definition}[theorem]{Definition}
\newtheorem{remark}[theorem]{Remark}

\numberwithin{equation}{theorem}

\allowdisplaybreaks[1]
\newcommand{\ds}{\displaystyle}

\DeclareMathOperator{\Pic}{Pic}

\newcommand\mm{\mathcal{M}}
\newcommand\mc{\mathcal{C}}
\newcommand\ml{\mathcal{L}}

\newcommand\C{\mathbb{C}}

\newcommand\SL{\operatorname{SL}}

\newcommand\frh{\mathfrak{h}}
\newcommand\frg{\mathfrak{g}}
\newcommand\frb{\mathfrak{b}}

\newcommand\frA{\mathfrak{A}}

\newcommand\frp{\mathfrak{p}}

\begin{document}

\title{On a lower bound for the dimension of non--abelian theta functions of positive genus}
\author{Arzu Boysal}
\date{}


\begin{abstract}

Let $\mc_g$ be a smooth projective irreducible curve over $\C$ of genus $g \geq 1$ and let
$\{p_1, p_2,\dots, p_s\}$ be a set of distinct points on $\mc_g$.  We fix a nonnegative
integer $\ell$ and denote by
$M_g(\underline{p},\underline{\lambda})$ the moduli space of parabolic semistable vector
bundles of rank $r$ with trivial determinant and fixed parabolic structure of
type $\underline{\lambda}=(\lambda_1,\lambda_2,\dots, \lambda_s)$ at $\underline{p}=(p_1, p_2,\dots, p_s)$,
where each weight $\lambda_i$ is in $P_{\ell}(\SL(r))$.
On $M_g(\underline{p},\underline{\lambda})$ there is a
canonical line bundle $\ml(\underline{\lambda}, \ell)$, whose sections are called
generalized parabolic $\SL(r)$-theta functions of order $\ell$.

The main result of this paper is: if $\sum_{1}^{s} \lambda_i$ is in the root lattice, then
$$\mbox{dim}\,H^0(M_g(\underline{p},\underline{\lambda}) , \ml(\underline{\lambda}, \ell))
\geq (\sharp P_{\ell}(\SL(r)))^{g-1},$$
where $\sharp P_{\ell}(\SL(r))=(\ell+r-1)!/(\ell ! (r-1)!)$.
Such nontrivial lower bounds for the
number of generalized parabolic theta functions were not previously known.

\end{abstract}

\maketitle

\pagestyle{myheadings}

\markboth{A.~Boysal}{Lower bound for the dimension of non--abelian theta functions}

\maketitle

\vspace{-1cm}
\section{Introduction}

Let $\mc_g$ be a smooth projective irreducible curve over $\C$ of genus $g\geq 1$
and let $\{p_1, p_2,\dots, p_s\}$ be a set of distinct points on $\mc_g$.
We fix a nonnegative integer $\ell$ and denote by
$M_g^{\SL(r)}(\underline{p},\underline{\lambda})$ (in short $M_g (\underline{p},\underline{\lambda})$)
the moduli space of parabolic semistable vector bundles of rank $r$ with trivial determinant and fixed parabolic
structure of type $\underline{\lambda}=(\lambda_1,\lambda_2,\dots, \lambda_s)$ at $\underline{p}=(p_1, p_2,\dots, p_s)$
(see [MS] or [LS] for definitions), where each weight $\lambda_i$ is in $P_{\ell}(\SL(r))$ (cf.~section 2.3).
For $g\geq 1$, $M_g(\underline{p},\underline{\lambda})$
is known to be nonempty for any such collection of weights $\underline{\lambda}$ (cf.~Remark~\ref{R:mnempty}).

On $M_g(\underline{p},\underline{\lambda})$ there is a
canonical line bundle $\ml(\underline{\lambda}, \ell)$, whose sections are called
generalized parabolic $\SL(r)$-theta functions of order $\ell$.  The dimension of
$H^0(M_g(\underline{p},\underline{\lambda}) , \ml(\underline{\lambda}, \ell))$ is given by the
celebrated Verlinde formula [V], but it does not give bounds, or even nonzeroness, due to signs involved.

In this paper we prove the following:

\begin{theorem}\label{T:main}
Let $\mc_g$ be a smooth projective irreducible curve over $\C$ of genus $g\geq 1$ and
let $\{p_1, p_2,\cdots, p_s\}$ be a set of distinct points on $\mc_g$.
Fix a nonnegative integer $\ell$ and denote by $M_g(\underline{p},\underline{\lambda})$
the moduli space of parabolic semistable $\SL(r)$-bundles
of type $\underline{\lambda}=(\lambda_1,\lambda_2,\dots, \lambda_s)$ at $\underline{p}=(p_1, p_2,\dots, p_s)$, with
each $\lambda_i$ in $P_{\ell}(\SL(r))$.
If $\sum_{1}^s{\lambda_i}$ lies in the root lattice, then
$$\mbox{dim}\, H^0(M_g(\underline{p},\underline{\lambda}), \ml(\underline{\lambda},\ell))\geq
(\sharp P_{\ell}(\SL(r)))^{g-1},$$
where $\sharp P_{\ell}(\SL(r))=(\ell+r-1)!/(\ell ! (r-1)!)$.
\end{theorem}

For $g=0$ and $M_0(\underline{p},\underline{\lambda})\neq \emptyset$,
dimension of $H^0(M_0(\underline{p},\underline{\lambda}), \ml(\underline{\lambda},\ell))$ is known to be
positive by works of Belkale on quantum horn and saturation conjectures [Be1].

In the proof of this main result we make essential use of the canonical identification of sections of
$\ml(\underline{\lambda}, \ell)$
with the space of conformal blocks (as given by Tsuchiya--Ueno--Yamada [TUY] via rational conformal field theory)
which is proved for the non-parabolic case (i.e.~s=0) independently by Beauville--Laszlo [BL], Faltings [F] and
Kumar--Narasimhan--Ramanathan [KNR]; parabolic variants are given by Pauly [P] and Laszlo--Sorger [LS].  For the
history of the problem and more detailed references we refer to Bourbaki article of Sorger [S].

Outline of the proof is as follows.  Under the above mentioned correspondence, the dimension
of  $H^0(M_g(\underline{p},\underline{\lambda}), \ml(\underline{\lambda},\ell))$, $n_g(\underline{\lambda})$,
obeys the `factorization rule' and has the `propagation of vacua' property as given by [TUY]
(cf.~Proposition~\ref{P:pvfr}).
Using these, we reduce the problem of determining nontrivial lower bounds for the dimension of sections  of
$\ml(\underline{\lambda},\ell)$ to the case of genus $g=1$ (cf.~Lemma~\ref{L:g1}).

For $g=1$, we first prove the existence of a nonzero section for the particular case of one marking (i.e.~$s=1$)
using a Parthasarathy-Ranga Rao-Varadaryan
(PRV) type result for the {\it fusion product} of $\SL(r)$ representations given by Belkale.  We show that
for any $\lambda \in P_{\ell}(\SL(r))$ that is also in the root lattice, there exists a
PRV weight $\mu \in P_{\ell}(\SL(r))$ such that, under
the degeneration of the curve to genus $0$, the canonical line with the `deformed' parabolic type
$(\lambda,\mu,\mu^*)$ admits a nonzero section (cf.~Proposition~\ref{P:onemark}).
We then reduce the arbitrary number of markings to this case (cf.~Proposition~\ref{P:torus})
by using again the PRV result as above and some properties of the fusion product.
By these two propositions,
we have $H^0(M_1(\underline{p},\underline{\lambda}), \ml(\underline{\lambda},\ell))\neq 0$ (that is
$n_1(\underline{\lambda})>0$, which is precisely the claim in Theorem~\ref{T:main} for $g=1$).
Then, by virtue of the above mentioned reduction to genus $1$, Theorem~\ref{T:main} follows.

The organization of this paper is as follows:  In section $2$, we set up the notation and give some preliminaries,
while in section $3$ we summarize basic definitions and results on the moduli stack and moduli
space of parabolic bundles, including the definition of the canonical line bundle $\ml(\underline{\lambda}, \ell)$.
In section $4$, we recall the relation between generalized parabolic theta functions and conformal blocks, and
state how the latter behave under degeneration of the curve.  In section $5$, we define the fusion product
and state the PRV result of Belkale.
Finally in section $6$, we prove the above mentioned reduction and propositions, and deduce the validity of the claim
in Theorem~\ref{T:main}.
\\
\\
\noindent {\bf Acknowledgements.} This work is in response to a
question of P.~Belkale, which asked if the dimension of
$H^0(M_g(\underline{p},\underline{\lambda}),
\ml(\underline{\lambda},\ell))$ is nonzero for $g\geq 1$.  I am very
thankful to him for many useful discussions, suggestions and reading
several earlier versions of this manuscript.  I also thank S.~Kumar,
M.S.~Narasimhan and T.R.~Ramadas.


\section{Preliminaries}

We set up the notation and recall basic definitions and concepts on simple algebraic groups.
Main reference for this section is Bourbaki [B].  The reader who is familiar with the notation
might as well skip to Section 3.

Let $G$ be a connected, simply connected, simple affine algebraic group over $\C$.
This will be our assumption on $G$ unless otherwise stated.

{\bf 2.1.}
We fix a Borel subgroup $B$ of $G$ and a maximal torus $T\subset B$. Let $\frh$,$\frb$ and $\frg$ denote the Lie
algebra of $T$, $B$ and $G$ respectively. Let $R=R(\frh,\frg)\subset \frh^*$ be the root system; there is the root
space decomposition $\frg=\frh \oplus (\oplus_{\alpha \in R}\frg_{\alpha})$.
We fix a basis $\Delta=\{\alpha_1, \dots, \alpha_r\}$ of $R$, where $r$ is the rank of $G$.  Let $R^+$ be the
set of positive roots.  Let $\frh_{\mathbb{R}}$ denote the real span of elements of $\frh$ dual to $\Delta$.
For each root $\alpha$, denote by $H_{\alpha}$ the unique element of
$[\frg_{\alpha},\frg_{-\alpha}]$ such that $\alpha(H_{\alpha})=2$; we denote by $\check{Q}$ the
lattice spanned by $H_{\alpha}$.
Let $\{\omega_i \}_{1\leq i\leq r}$  be the set of fundamental weights, defined as the basis of $\frh^*$
dual to $\{H_{\alpha_i}\}_{1\leq i\leq r}$.  Denote the weight lattice by $P \subset \frh^*$, that is,
$P=\{ \lambda \in \frh^* : \lambda(H_{\alpha}) \in \mathbb{Z},\; \forall \alpha \in R \}$.

{\bf 2.2.}
Let $(\;|\;)$ denote the Killing form on $\frg$ normalized such that $(H_{\theta}|H_{\theta})=2$, where
$\theta$ is the highest root of $G$.  We will use the same notation for the restricted form on $\frh$, and
the induced form on $\frh^*$.
We introduce the Weyl group $W:=N_G(T)/T$.  We fix a (closed) positive Weyl chamber
in $\frh$, $\frh_+:=\{x\in \frh_{\mathbb{R}}: \alpha(x)\geq 0\;\forall \alpha \in R^+\}$;
correspondingly a positive Weyl chamber in $\frh^*$,
$\frh_{+}^*:=\{\omega \in \frh_{\mathbb{R}}^*  : \omega(H_{\alpha})\geq 0 \;\forall \alpha \in R^+\}.$
The Killing form is positive definite on $\frh_{\mathbb{R}}$, thus induces the following canonical isomorphism
\[F: \frh_{\mathbb{R}}^* \to \frh_{\mathbb{R}}\;;\;\; \alpha \mapsto (2/(H_{\alpha}| H_{\alpha}))H_{\alpha}\]
with $F^{-1}(H_\alpha)=(2/(\alpha|\alpha))\alpha$.
Under the above identification, $W$ is generated by elements $w_{\alpha}$ for $\alpha \in \Delta$; acting
on $\frh_{\mathbb{R}}^*$ as $w_{\alpha}\beta=\beta-\beta(H_{\alpha})\alpha$.
Thus, given a weight $\beta \in P$ that is also in the root lattice, for any $w \in W$, $w\beta-\beta$ is also in the
root lattice.
Let $w_0$ denote the
longest element in the Weyl group.  Then for any $\lambda \in P$, the dual $\lambda^*=-w_0\lambda$.

{\bf 2.3.}
Linear combinations of $\{\omega_i\}_{1\leq i\leq r}$ that are in $\frh_{+}^*$ with integer coefficients are called
dominant weights; we denote the set of dominant weights by $P^+$.
For a positive integer $\ell$, define
\[P_{\ell}(G):=\{\lambda \in P^+ : \lambda(H_{\theta})\leq \ell\}.\]
The set $P_{\ell}(G)$ is finite and stable under taking the dual.  In particular, for
$G=\SL(r)$ and $\lambda=\sum_{i=1}^{r-1}a_i\omega_i$,
\[P_{\ell}(\SL(r))=\{a_i \in \mathbb{Z}_{\geq 0} : \sum_{i=1}^{r-1} a_i\leq \ell\}.\]

{\bf 2.4.}
Let \[\mathfrak{A}=\{h\in \frh_{+} : \theta(h)\leq 1\}\] denote the (closed) fundamental alcove.
$\frA$ is a fundamental domain for the action on $\frh_{\mathbb{R}}$ of the affine Weyl
group $W_{\mbox{\scriptsize{aff}}}$, which is the semidirect product of  $W$ and $\check{Q}$.
By virtue of the normalization in the Killing form, for any $\lambda \in P_{\ell}(G)$,
\[\theta(F(\lambda))=\lambda(F(\theta))=\lambda(H_{\theta}).\]
Thus, $\lambda/\ell$ determines a unique element $F(\lambda/\ell)$ in $\frA$.
Often we shall implicity identify a weight in $P_{\ell}(G)$ with the corresponding element in $\frA$, and hence
with an orbit of $\frh_{\mathbb{R}}/W_{\mbox{\scriptsize{aff}}}$.

{\bf 2.5.}
A subset $S$ of $\Delta$ defines a subalgebra
$\frp_{S}=\frb \oplus (\oplus_{\alpha \in S}\frg_{-\alpha})\subset \frg$, and hence a
standard parabolic subgroup $P_S \subset G$.  A dominant weight $\lambda$ of $G$ uniquely determines a
standard parabolic subgroup $P\subset G$ and a dominant character $\chi \in X(P)$, where $X(P)$ denotes the
character group of $P$.


\section{Review of definitions and results on moduli stack and moduli space of parabolic bundles
over $\mc_g$}

Main reference for this section is Laszlo--Sorger [LS]; also see Mehta--Seshadri
[MS] for the case $G=\SL(r)$.

Let $\mc_g$ be a smooth projective irreducible curve over $\C$ of genus $g \geq 1$ and
let $\underline{p}=(p_i)_{1\leq i \leq s}$ be distinct points on $\mc_g$ labelled by
dominant weights $\underline{\lambda}=(\lambda_i)_{1\leq i \leq s}$.
We fix a base point $b$ on $\mc_g$ distinct from $\underline{p}$. For each weight $\lambda_i$
define a subset $S_i$ of $\Delta$ as $S_i:=\{\alpha \in \Delta:\lambda_i(H_{\alpha})=0\}$.
We denote by $P_i$ the standard parabolic subgroup associated to $S_i$ (as in section 2.5).
For any $i$, $1\leq i \leq s$, and $j$ such that $\alpha_j \in \Delta \diagdown S_i$,
let $t_{i,j}=\lambda_i(H_{\alpha_j})$ denote the corresponding parabolic weight.

\begin{definition}
(i). A {\it quasi-parabolic $G$-bundle} of type $\underline{P}=(P_1,\dots, P_s)$ at $\underline{p}=(p_1,\dots,p_s)$ is a
$G$-bundle $E$ on $\mc_g$ together with a fibre $F_i$ of the associated bundle $E(G/P_i)$ at $p_i$ for all $i$,
$1\leq i \leq s$.

(ii). A {\it parabolic $G$-bundle} of type $(\underline{P},\underline{t})$ at $\underline{p}$ is a
quasi-parabolic $G$-bundle of type $\underline{P}$ at $\underline{p}$, together with parabolic weights
$\{t_{i,j}\}$ for each $i$, $1 \leq i \leq s$, and $j$ such that $\alpha_j \in \Delta \diagdown S_i$.

(iii). A {\it family of quasi-parabolic G-bundles} of type $\underline{P}$ on $\mc_g$,
parametrized by a $\mathbb{C}$-scheme $Z$, denoted by $(E,\underline{\sigma})$, is a $G$-bundle $E$ over
$Z \times \mc_g$ together
with sections $\sigma_i:Z\to E(G/P_i)|_{Z \times \{p_i\}}$ for each $i$, $1\leq i \leq s$.

(iv). A {\it morphism} from $(E,\underline{\sigma})$ to $(E',\underline{\sigma'})$ is a morphism
$f:E\to E'$ of G-bundles such that for any $i$, $1 \leq i \leq s$, the restriction of $f$ to the fibre over
$Z\times \{p_i\}$, $f|_{Z\times \{p_i\}}$,
satisfies $\underline{\sigma'}=f|_{Z\times \{p_i\}}\circ \underline{\sigma}$.

\end{definition}

Let $\mm_{g} ^G(\underline{p},\underline{\lambda})$ denote the (algebraic) moduli stack of
quasi-parabolic $G$-bundles of type $\underline{P}$ at $\underline{p}$ determined by $\underline{\lambda}$.
Its objects are families of quasiparabolic G-bundles of type $\underline{P}$ parametrized by $Z$, and morphisms
are isomorphisms between such families.


It is known that (though we won't strictly need this, it helps to clarify the form of the
canonical line bundle that we will shortly define),
\[\Pic(\mm_g ^G(\underline{p},\underline{\lambda}))=
\mathbb{Z}\mathcal{L}\times \prod_{i=1}^s X(P_i),\]
where for the particular case $G=\SL(r)$, $\mathcal{L}$ is the determinant line bundle.
[LS, Theorem 1.1].

For a family $(E,\underline{\sigma})$ of type $\underline {P}$ parametrized by $Z$, and
for fixed $i$ and $j$ ($1 \leq i \leq s$ and $j$ such that $\alpha_j \in \Delta \diagdown S_i$)
one can define a line bundle on $Z$ as the pullback of the
associated line bundle $E\times_{P_i}\mathbb{C}_{-\omega_j}$ on $E(G/P_i)$ (see $E \to E(G/P_i)$
as a $P_i$-bundle) using the section $\sigma_i$.  Let $\mathcal{L}_{i,j}$ denote the line bundle on the
stack, constructed as such on any $Z$.  Then, there is a canonical line bundle on
$\mm_g ^G(\underline{p},\underline{\lambda})$ defined as
\[\mathcal{L}(\underline{\lambda},\ell)=\mathcal{L}^{\ell}\boxtimes
(\boxtimes_{i=1}^n(\boxtimes_{j :\; \alpha_j\in \Delta \diagdown S_i}
\mathcal{L}_{i,j}^{t_{i,j}}))\] where $\{t_{i,j}\}$ are the parabolic weights (which are positive
integers by definition).

There exists a (coarse) moduli space for $\mm_g ^G(\underline{p},\underline{\lambda})$, which is
a projective variety and its points are equivalence classes of parabolic semistable G-bundles
of a fixed topological type and fixed parabolic structure [BR, Theorem II].
We denote this space by $M_g ^G(\underline{p},\underline{\lambda})$.

Fix a maximal compact subgroup $K$ of $G$.  Let $C_g$ denote the Riemann surface with genus $g$.
The real analytic space underlying $M_g ^G(\underline{p},\underline{\lambda})$
admits a description as the space of representations of the fundamental group
$\pi_1(C_g \setminus \{p_1,\dots, p_s\},b)$ into $K$ upto conjugation.
(For the proof of this identification we refer to Bhosle--Ramanathan [BR, Proposition 2.3], see also Teleman--Woodward [TW];
the particular case of $G=\SL(r)$ is given by Mehta--Seshadri [MS, Theorem 4.1].
In [MS], the case $K=U(r)$ and $g\geq 2$ is considered, but the generalization of the above identification to
$K=SU(r)$ and $g \geq 0$ is immediate and well known).

Recall that $\pi_1(C_g \setminus \{p_1,\dots, p_s\},b)$ is a free group on $2g+s$ generators
with one relation:
\[ \begin{aligned}
\pi_1(C_g \setminus \{p_1,\dots, p_s\},b)=&<a_1,a_2,\dots,a_g,b_1,b_2,\dots,b_g,c_1,\dots,c_s: \\
&\prod_{i=1}^g [a_i,b_i]=\prod_{j=1}^s c_j>,\end{aligned} \]
where $[a,b]:=aba^{-1}b^{-1}$ for any symbol $a$ and $b$.
As is well known, the set of conjugacy classes of $K$ is in bijective correspondence with
$\frh_{\mathbb{R}}/W_{\mbox{\scriptsize{aff}}}$.
For a collection of weights $\underline{\lambda}=(\lambda_1,\dots,\lambda_s)$ in $P_{\ell}(G)$, we denote
by $C_i$ the conjugacy class of $\mbox{exp}(2 \pi \sqrt{-1}F(\lambda_i/\ell))$, $1\leq i \leq s$
(see section 2.2 for the definition of $F$ and 2.4 for the assignment).
Then, we have the following isomorphism as real analytic spaces
\[
\begin{aligned}
M_g ^G(\underline{p},\underline{\lambda})\simeq &\{(k_1,\dots, k_{2g}, h_1,\dots,h_s) \in K^{2g+s}:\\
&\prod_{i=1}^g [k_i,k_{i+g}]=\prod_{j=1}^s h_j,\;\;h_i \in C_i\}/\mbox{Ad}K \;\;\;\;\;\;(*),
\end{aligned} \]
where $/\mbox{Ad}K$ refers to the quotient under the diagonal adjoint action of $K$ on $K^{2g+s}$.

\begin{remark}\label{R:mnempty}By a theorem of Got\^{o} [G], the image of the commutator map
$K \times K \to K  ; \;(k_1,k_2)\mapsto [k_1,k_2]$ is surjective.
Thus, under the identification $(*)$ above, for $g\geq 1$ the moduli space $M_g ^G(\underline{p},\underline{\lambda})$ is
nonempty for any collection of weights $\underline{\lambda}$ in $P_{\ell}(G)$.
\end{remark}

\begin{remark}\label{R:descent}
For $G=\SL(r)$, $\mathcal{L}(\underline{\lambda},\ell)$ descends to
$M_g ^G(\underline{p},\underline{\lambda})$ (see for example [DN] or [P]).
\end{remark}

Moreover, the following canonical correspondence exists between their sections.
\begin{theorem}\label{T:section}
$$H^0(\mm_g ^{\SL(r)}(\underline{p},\underline{\lambda}),\mathcal{L}(\underline{\lambda},\ell) )\simeq
H^0(M_g ^{\SL(r)}(\underline{p},\underline{\lambda}), \mathcal{L}(\underline{\lambda},\ell)).$$
\end{theorem}

For the non-parabolic case the proof of Theorem~\ref{T:section} can be found in Beauville--Laszlo [BL, Propositions 8.3
and 8.4]; the
generalization to the parabolic case is immediate and given by Pauly [P, Proposition 5.2].  For any $G$ simple
simply connected see Kumar--Narasimhan--Ramanathan [KNR, Propositions 6.4 and 6.5, and Theorem 6.6] and
Narasimhan--Ramadas [NR].


\section{Generalized parabolic theta functions and conformal blocks}

In this section we will recall the relation between generalized parabolic $G$-theta functions of level $\ell$,
that is,  $H^0(\mm_g ^G(\underline{p},\underline{\lambda}), \mathcal{L}(\underline{\lambda},\ell))$, and
the space of conformal blocks given by Tsuchiya, Ueno and Yamada [TUY].

Let $\tilde{\frg}=\frg \otimes \mathbb{C}((z))\oplus \mathbb{C}\mbox{C}$
denote the (untwisted) affine lie algebra associated to $\frg$ over $\mathbb{C}((z))$, with the lie bracket
given by $$[x\otimes f,y\otimes g]=[x,y]\otimes fg+((x|y)Res(gdf))\cdot \mbox{C}\;\mbox{and}\;[\frg,\mbox{C}]=0$$ for
$x,y \in \frg$ and $f,g \in \mathbb{C}((z))$.

We fix a positive integer $\ell$, called the {\it level}, and recall the definition of
$P_{\ell}(G)$ from section 2.3.
For each $\lambda \in P_{\ell}(G)$ there exists a unique (upto isomorphism) left $\tilde{\frg}$-module
$\mathcal{H}_{\lambda}$, called the integrable highest weight $\tilde{\frg}$-module, with the central
element $\mbox{C}$ acting as $\ell\cdot \mbox{Id}$ (see Kac [K] for further details).

Let $\mathfrak{a}$ be any lie algebra; denote by $\mathcal{U}(\mathfrak{a})$ its universal enveloping algebra and by
$\mathcal{U^+}(\mathfrak{a}):=\mathfrak{a}\mathcal{U}(\mathfrak{a})$ the augmentation ideal of $\mathcal{U}(\mathfrak{a})$.
For an $\mathfrak{a}$-module $V$, the quotient module $V/\mathcal{U^+}(\mathfrak{a})V$ is said to be
the {\it space of coinvariants} of $V$ with respect to $\mathfrak{a}$, and is denoted by $[V]_{\mathfrak{a}}$.

\begin{definition}\label{D:cb}
The space of {\it conformal blocks} on $\mc_g$ with marked points $(p_1,\dots,p_k)$ and weights
$(\lambda_1,\cdots,\lambda_k)$ attached to them
(each $\lambda_i \in P_{\ell}(G)$) with central charge $\ell$ is  $$V_{\mc_g} ^G(\underline{p},\underline{\lambda}):=
[\mathcal{H}_{\lambda_1}\otimes \mathcal{H}_{\lambda_2}\otimes\cdots\otimes \mathcal{H}_{\lambda_k}]
_{\frg\otimes\mathcal{O}(\mc_g\backslash\{p_1,\dots,p_k\})},$$ where $\mathcal{O}(\mc_g\backslash\{p_1,\dots,p_k\})$
denotes the ring of algebraic functions on the punctured curve $\mc_g\backslash\{p_1,\dots,p_k\}$.
(See [TUY] for details of the particular action used in taking the coinvariants.)
\end{definition}

These vector spaces form a projectively flat vector bundle over the moduli space of smooth projective complex curves
with marked points;
hence the rank of the vector space $V_{\mc_g} ^G(\underline{p},\underline{\lambda})$ depends only on the genus
$g$ and on weights $\underline{\lambda}$ [TUY, Lemma 2.3.2 and Remark 4.1.7].  Denote this common rank by
$n_g(\underline{\lambda})$, that is, $n_g(\underline{\lambda})=\mbox{dim}V_{\mc_g} ^G(\underline{p},\underline{\lambda})$.
Its value is given by the Verlinde formula [V].

Below, we state the principle of `propagation of vacua' and the `factorization rule' for conformal blocks,
only in terms of their dimensions, as we won't need explicit isomorphisms.

\begin{proposition}\label{P:pvfr}
For any collection of weights $\underline{\lambda}$ in $P_{\ell}(G)$:
\[\begin{array}{l}
(a). (\mbox{Propogation of vacua, [TUY, Proposition 2.2.3]})\;\; n_{g}(\underline{\lambda},0)=n_{g}(\underline{\lambda}).\\
(b). (\mbox{Factorization rule, [TUY, Proposition 2.2.6]})\\
    \;\;\;\;\;\; n_g(\underline{\lambda})=\ds \sum_{\nu \in P_{\ell}(G)}n_{g-1}(\underline{\lambda},\nu, \nu^*).
\end{array}\]
\end{proposition}

\begin{theorem}\label{T:confb}
There is a canonical isomorphism
$$H^0(\mm_g ^G(\underline{p},\underline{\lambda}), \mathcal{L}(\underline{\lambda},\ell))\simeq
V_{\mc_g} ^G(\underline{p},\underline{\lambda}).$$
\end{theorem}

See introduction for the list of references attributed to the proof of Theorem~\ref{T:confb}.

Finally, by Remark~\ref{R:descent}, Theorem~\ref{T:section} and Theorem~\ref{T:confb}, for $G=\SL(r)$,

$$H^0(M_g ^{\SL(r)}(\underline{p},\underline{\lambda}), \mathcal{L}(\underline{\lambda},\ell))\simeq
V_{\mc_g} ^{\SL(r)}(\underline{p},\underline{\lambda})\;\;(**).$$


\section{Fusion product and PRV components}

Let $\mathcal{R}(\frg)$ denote the Grothendieck ring of finite dimensional representations of $\frg$.
There is a bijection of $P^+$ onto $\mathcal{R}(\frg)$,
where a dominant weight $\lambda$ is associated to the isomorphism class of the simple $\frg$-module $V_{\lambda}$
containing a highest weight vector with weight $\lambda$.  The multiplicative structure in $\mathcal{R}(\frg)$
is induced from the tensor product of two representations.

The {\it fusion ring} associated to $\mathfrak{g}$ and
the nonnegative integer $\ell$, $\mathcal{R}_{\ell}(\mathfrak{g})$,
is a free $\mathbb{Z}$-module with basis $\{ V_{\lambda}, \lambda \in P_{\ell}(G) \}$.
The ring structure is defined as follows:

$$V_{\lambda}\otimes ^F V_{\mu}:= \ds \bigoplus_{\nu \in P_{\ell}(G)} n_{\lambda,\mu}(\nu^*) V_{\nu},$$
where
$$n_{\lambda,\mu}(\nu^*):=\mbox{dim}[ \mathcal{H}_{\lambda}\otimes \mathcal{H}_{\mu} \otimes \mathcal{H}_{\nu} ^* ]
_{\mathfrak{g}\otimes \mathcal{O}(\mathbb{P}^1 - \{\mbox{3 points} \})}=n_0(\lambda,\mu,\nu^*),$$
where $\mathcal{H}_{\lambda}$ refers to the irreducible integrable $\tilde{\frg}$-module of highest weight
$\lambda$ as described in section $4$.  The product $\otimes^F$ is clearly commutative; it is also
associative by virtue of the `factorization rule'.

We now recall (a version of) the Parthasarathy-Ranga Rao-Varadaryan (PRV)
conjecture (already proven) for the tensor product of two simple $\frg$-modules.
Let $V_{\lambda_1}$ and $V_{\lambda_2}$ be two finite dimensional irreducible $\mathfrak{g}$-modules
with highest weights $\lambda_1$ and $\lambda_2$ respectively. Then, for any $w \in W$
the irreducible $\mathfrak{g}$-module $V_{\overline{\lambda_1+w \lambda_2}}$  occurs with multiplicity at
least one in $V_{\lambda_1}\otimes V_{\lambda_2}$, where $\overline{\lambda_1+w \lambda_2}$ denotes the unique
dominant element in the $W-$orbit of $\lambda_1+w \lambda_2$.  We refer to Kumar [Ku] for the proof.

An analogous theorem for the {\it fusion product} of $\mathfrak{sl}(r)$ representations is given by Belkale.
Before stating this theorem we observe the following.

Given any $\lambda_1, \lambda_2 \in P_{\ell}(G)$ and $w \in W_{\mbox{\scriptsize{aff}}}$,
there exists $\tilde{w} \in W_{\mbox{\scriptsize{aff}}}$ such that $\tilde{w}(F(\lambda_1/\ell)+
w F(\lambda_2/\ell))\in \frA$ is the unique element in the $W_{\mbox{\scriptsize{aff}}}-$orbit of
$F(\lambda_1/\ell)+w F(\lambda_2/\ell)$ (see section 2.4).
Clearly,
$$F^{-1}(\tilde{w}(F(\lambda_1/\ell)+wF(\lambda_2/\ell)))=\tilde{w}(\lambda_1/\ell +w\lambda_2/\ell).$$
Moreover, $\ell(\tilde{w}(\lambda_1/\ell +w\lambda_2/\ell))$ is in $P_{\ell}(G)$ by construction, and is of the form
$\tilde{w}'(\lambda_1 +w\lambda_2)$ for some $\tilde{w}' \in W_{\mbox{\scriptsize{aff}}}$ such that
$\tilde{w}' \equiv \tilde{w} \;\;(\mbox{mod}\;\ell\check{Q})$.
We denote this element in $P_{\ell}(G)$ (uniquely determined by the triple $\lambda_1$, $\lambda_2$ and $w$)
by $\overline{\lambda_1+w \lambda_2}^F$.

\begin{theorem}([Be1], Theorem 1.7)~\label{T:gzero}  With notation as before,
$M_0^{\SL(r)}(\underline{p},\underline{\lambda})\neq \emptyset$ if and only if
$n_0(\underline{\lambda})\neq 0$.
\end{theorem}

\begin{theorem}\label{T:prv} (PRV for fusion product, [Be2])
Let $V_{\lambda_1}$ and $V_{\lambda_2}$ be finite dimensional irreducible $\mathfrak{sl}(r)$-modules with
highest weights $\lambda_1$ and $\lambda_2$ respectively, with $\lambda_1,\lambda_2 \in P_{\ell}(\SL(r))$.
Then, for any $w \in W_{\mbox{\scriptsize{aff}}}$, the irreducible $\mathfrak{sl}(r)$-module
$V_{\overline{{\lambda_1+w\lambda_2}}^F}$ occurs with multiplicity at least one in
$V_{\lambda_1}\otimes ^F V_{\lambda_2}$.
\end{theorem}

\begin{proof}
(P. Belkale informed me that he learned of this derivation of PRV from Theorem 1.7. in [Be1]
from C. Woodward).

Let $\mu_i=F(\lambda_i/\ell)$ for $i=1,2$, and
$\mu_{3}=\tilde{w}F(\lambda_1/\ell+w\lambda_2/\ell)=\tilde{w}(\mu_1+w\mu_2)$, where
$\tilde{w}(\mu_1+w\mu_2) \in \frA$ is the unique element in the $W_{\mbox{\scriptsize{aff}}}$-orbit
of $\mu_1+w\mu_2$.
Clearly, there exist elements $w_1,w_2, \; \mbox{and} \;w_3 \in W_{\mbox{\scriptsize{aff}}}$ such that
$w_1\mu_1+ w_2\mu_2 +w_3\mu_3^*=0$, where $\mu_3^*$ denotes the dual of $\mu_3$ as in section 2.2.
Thus, under the identification $(*)$ in section 3,
$M_0 ^{\SL(r)}(\underline{p},\lambda_1,\lambda_2,(\overline{{\lambda_1+w \lambda_2}}^F)^*)\neq \emptyset$.
Then it follows from Theorem~\ref{T:gzero} that $n_0(\lambda_1,\lambda_2,(\overline{{\lambda_1+w\lambda_2}}^F)^*)\neq 0$;
equivalently $V_{\overline{{\lambda_1+w\lambda_2}}^F} \subset V_{\lambda_1}\otimes ^F V_{\lambda_2}$ by the definition of
the fusion product.
\end{proof}


\section{Proof of Theorem~\ref{T:main}}

In this section we restrict ourselves to  $G=\SL(r)$.
Recall that under the identification $(**)$ of section 4,
\[\mbox{dim}\,H^0(M_g ^{\SL(r)}(\underline{p},\underline{\lambda}), \mathcal{L}(\underline{\lambda},\ell))=
\mbox{dim}\,V_{\mc_g} ^{\SL(r)}(\underline{p},\underline{\lambda})=n_g(\underline{\lambda}) \tag{1}.\]

\begin{lemma}\label{L:g1}
It suffices to show that Theorem~\ref{T:main} holds for $g=1$.
\end{lemma}
\begin{proof}
Using the `factorization rule', Proposition~\ref{P:pvfr}$(b)$, for any
collection of weights $\underline{\lambda}$ in $P_{\ell}(\SL(r))$
\[ \begin{array}{ll}
n_g(\underline{\lambda})&=\ds \sum_{\nu \in P_{\ell}(\SL(r))}n_{g-1}(\underline{\lambda},\nu, \nu^*) \\
&=n_{g-1}(\underline{\lambda},0,0)+\ds \sum_{\nu \in P_{\ell}(\SL(r)),\;\nu\neq 0}n_{g-1}(\underline{\lambda},\nu, \nu^*).
\end{array}\]
By `propagation of vacua', Proposition~\ref{P:pvfr}$(a)$,
$n_{g-1}(\underline{\lambda},0,0)=n_{g-1}(\underline{\lambda})$.  Hence we get \[ n_{g}(\underline{\lambda})\geq
n_{g-1}(\underline{\lambda})\tag{2}.\]

Now suppose $n_{1}(\underline{\lambda})>0$. Clearly, by inequality $(2)$,
$n_{g}(\underline{\lambda})>0$ for any $g\geq 1$.  Moreover, since
$\nu+\nu^*$ is in the root lattice for any $\nu \in P$,
we get the lower bound as claimed in Theorem~\ref{T:main} (under identification $(1)$) inductively by
Proposition~\ref{P:pvfr}$(b)$.
\end{proof}

\begin{remark} In the non-parabolic case, that is when $s=0$,
$n_1(\emptyset)$ is equal to the dimension of the fusion ring (i.e.~ cardinality of $P_{\ell}(\SL(r))$
by the Verlinde formula [V]) which is positive for any $\ell \geq 0$.  Thus, by virtue of Lemma~\ref{L:g1} and
equality $(1)$, Theorem~\ref{T:main} trivially holds for this particular case
(in fact for any $G$ under the assumptions of section 2 with the corresponding canonical line that descends).
\end{remark}

\begin{proposition}\label{P:onemark}
Theorem~\ref{T:main} holds for the case of one marking on $\mc_1$ (i.e. for the particular case $g=1$ and $s=1$).
\end{proposition}
\begin{proof}
We will show that for any $\lambda \in P_{\ell}(\SL(r))$ that is also in the root lattice
there exists $\mu \in P_{\ell}(\SL(r))$ such that $V_{\mu} \subset V_{\lambda} \otimes^F V_{\mu}$.
This will imply by the `factorization rule' that $n_1(\lambda)>0$.

Given $\lambda=\sum_{i=1}^{r-1} a_i \omega_i$ that lies in the root lattice, by change of basis
$\lambda=\sum_{i=1}^{r-1} n_i \alpha_i$ for some  $n_i \in \mathbb{Z}$.
For any $w\in W_{\mbox{\scriptsize{aff}}}$, consider the equation \[\lambda=\mu'-w \mu' \tag{3}.\]
There exits a $\mu' \in P$ satisfying $(3)$ if and only if $\lambda$ lies in the root lattice.
The latter granted, we choose $w=w_{\alpha_{r-1}}w_{\alpha_{r-2}}\cdots w_{\alpha_2}w_{\alpha_1} \in W$.
It is trivial to see that, for this particular choice of $w$, $\mu'=\sum_{i=1}^{r-1}(n_{i}-n_{i-1})\omega_i$
(with $n_0=0$) is the solution in $P$.
Let $\tilde{w} \in W_{\mbox{\scriptsize{aff}}}$ such that $\mu:=\tilde{w}^{-1}\mu'$ is in
$P_{\ell}(\SL(r))$.  Then we have
\[\lambda=\tilde{w}\mu-w\tilde{w}\mu,\]
that is, $\overline{\lambda+w \tilde{w} \mu}^F=\mu$.
Therefore, by Theorem~\ref{T:prv}, $V_{\mu} \subset V_{\lambda}\otimes ^F V_{\mu}$;
equivalently $n_0(\lambda,\mu,\mu^*)>0$.  Then it follows from Proposition~\ref{P:pvfr}$(b)$ that
$n_1(\lambda)>0$.  Hence the claim by equality $(1)$.
\end{proof}

\begin{proposition}\label{P:torus}
Theorem~\ref{T:main} holds for arbitrary number of markings on $\mathcal{C}_1$.
\end{proposition}

\begin{proof}
We will show that $n_{1}(\underline{\lambda})>0$ whenever $\sum_{i=1}^s \lambda_i$
lies in the root lattice of $\SL(r)$.

Recall that for any $\lambda_1, \lambda_2 \in P_{\ell}(\SL(r))$ and $w \in W_{\mbox{\scriptsize{aff}}}$,
there exists $\tilde{w}' \in W_{\mbox{\scriptsize{aff}}}$ such that
$\tilde{w}'(\lambda_1+w \lambda_2)$ is in $P_{\ell}(\SL(r))$.  Moreover, by Theorem~\ref{T:prv},
the $\mathfrak{sl}(r)$ representation
with highest weight $\tilde{w}'(\lambda_1+w \lambda_2)$ appears in $V(\lambda_1) \otimes^F V(\lambda_2)$.
Clearly,
$\tilde{w}'(\lambda_1+w\lambda_2)=\lambda_1+\lambda_2+(\tilde{w}'\lambda_1-\lambda_1)+
(\tilde{w}' w\lambda_2-\lambda_2).$
Thus, if $\lambda_1+\lambda_2$ is in the root lattice, so is $\tilde{w}'(\lambda_1+w\lambda_2)$ (see section 2.2).

We now repeat the above construction.  Given weights
$\lambda_1,\dots,\lambda_s \in P_{\ell}(\SL(r))$ such that $\sum_{i=1}^s \lambda_i$ is in the root lattice,
and $w_1 \in W_{\mbox{\scriptsize{aff}}}$, there exists
$w_2,w_3,\dots,w_s$ in $W_{\mbox{\scriptsize{aff}}}$ such that the weight
\[\tilde \lambda=w_s( \lambda_s+w_{s-1}(\lambda_{s-1}+w_{s-2}(\lambda_{s-2}+\cdots)))\]
is in $P_{\ell}(\SL(r))$ (and in the root lattice).  Moreover, by repeated application of Theorem~\ref{T:prv},
\[V_{\lambda_1}\otimes^F V_{\lambda_2}\otimes^F \cdots \otimes^F V_{\lambda_s}\supset V_{\tilde{\lambda}}.\tag{4}\]

Since $V_{\lambda_1}\otimes^F V_{\lambda_2}\otimes^F \cdots \otimes^F V_{\lambda_s}$ decomposes with positive
structure coefficients and the fusion product is associative (see section 5), for any $\nu \in P_{\ell}(\SL(r))$
\[(V_{\lambda_1}\otimes^F \cdots \otimes^F V_{\lambda_s})\otimes^F V_{\nu}=V_{\lambda_1}\otimes^F \cdots \otimes^F V_{\lambda_s}\otimes^F V_{\nu}
\supset V_{\tilde \lambda}\otimes^F V_{\nu}.\]
In particular,
\[n_0(\underline{\lambda}, \nu, \nu^*) \geq n_0(\tilde \lambda, \nu, \nu^*)\tag{5}\]
where $\underline{\lambda}=(\lambda_1,\dots,\lambda_s)$.

Using Proposition~\ref{P:pvfr}$(b)$,
\[ n_1(\underline{\lambda})=\sum_{\nu \in P_{\ell}(\SL(r))}n_0(\underline{\lambda}, \nu, \nu^*)\geq
\sum_{\nu \in P_{\ell}(\SL(r))}n_0(\tilde \lambda, \nu, \nu^*),\tag{6}\]
where the inequality in $(6)$ holds by $(5)$ term by term. Using one more time Proposition~\ref{P:pvfr}$(b)$,
\[\sum_{\nu \in P_{\ell}(\SL(r))}n_0(\tilde \lambda , \nu, \nu^*) =n_1(\tilde \lambda).\tag{7}\]
By construction $\tilde \lambda$ is in root lattice and in $P_{\ell}(\SL(r))$.
Thus, by Proposition~\ref{P:onemark}, $n_1(\tilde \lambda)>0$.  Then it follows from $(6)$ and $(7)$ that
$n_1(\underline{\lambda})>0$.  Hence the claim, under the identification $(1)$.
\end{proof}

By the virtue of Lemma~\ref{L:g1} and Proposition~\ref{P:torus}, Theorem~\ref{T:main} holds.


\section*{References}
\vspace{-0.5cm}
\[ \begin{array}{ll}

\mbox{[B]} & \mbox{Bourbaki N., {\it Groupes et Alg\`ebres de Lie}, Chap. 4--6, Masson, Paris,1981.}\\

\mbox{[Be1]} & \mbox{Belkale P., {\it Quantum generalization of the Horn conjecture}, math.AG/0303013.}\\
             & \mbox{latest version, http://www.unc.edu/~belkale/}\\

\mbox{[Be2]} &\mbox{Belkale P., {\it Private communication.}}\\

\mbox{[BL]}& \mbox{Beauville A. and Laszlo Y., {\it Conformal blocks and generalized theta functions},}\\
            & \mbox{Commun. Math. Phys. {\bf 164}, 385--419, (1994).}\\

\mbox{[BR]} &\mbox{Bhosle U. and Ramanathan A., {\it Moduli of Parabolic G-bundles on curves},}\\
           & \mbox{Math. Z. {\bf 202}, 161--180, (1989).}\\

\mbox{[DN]} & \mbox{Drezet J.--M. and Narasimhan M.S., {\it Groupe de Picard des vari\'et\'es de}}\\
            & \mbox{{\it modules de fibr\'es semi--stables sur les courbes alg\'ebriques}, Invent. Math.}\\
           & \mbox{{\bf 97}, 53--94, (1989).}\\

\mbox{[F]}& \mbox{Faltings G., {\it A proof for the Verlinde formula}, J. Alg. Geom. {\bf 3}, 347--374,}\\
           & \mbox{(1994).}\\

\mbox{[G]} &\mbox{Got\^{o} M., {\it A theorem on compact semi-simple groups}, J. Math. Soc. Japan {\bf 1},}\\
           & \mbox{no.3, 270--272, (1949).}\\

\mbox{[K]} &\mbox{Kac V., {\it Infinite dimensional lie algebras}, Cambridge University Press, }\\
            &\mbox{Cambridge}, 1990.\\

\mbox{[Ku]} & \mbox{Kumar S., {\it Proof of the Parthasarathy--Rango Rao--Varadarajan conjecture},}\\
            &\mbox{Invent. math. {\bf 93}, 117--130, (1988).}\\

\mbox{[KNR]} & \mbox{Kumar S., Narasimhan M.S. and Ramanathan A., {\it Infinite Grassmannians}}\\
            &\mbox{{\it and moduli spaces of $G$--bundles}, Math. Annalen {\bf 300}, 41--75, (1994).}\\

\mbox{[LS]} & \mbox{Laszlo Y. and Sorger C., {\it The line bundles on the moduli of parabolic}}\\
        &\mbox{{\it G--bundles over curves and their sections}, Ann. Scient. Ec. Norm. Sup.}\\
        &\mbox{$4^{e}$ s\`erie, t.30, 499--525, (1997).}\\

\mbox{[MS]} & \mbox{Mehta V.B. and Seshadri C.S., {\it Moduli of vector bundles on curves with}}\\
        &\mbox{{\it parabolic structures}, Math. Ann. {\bf 248}, 205--239, (1980).}\\

\mbox{[NR]} & \mbox{Narasimhan M.S. and Ramadas T.R., {\it Factorisation of generalised theta}}\\
      &\mbox{functions. I., Invent. Math. {\bf 114}, no.3, 565--623, (1993).}\\

\mbox{[P]} & \mbox{Pauly C., {\it Espaces de modules paraboliques et blocks conformes}, Duke}\\
        &\mbox{Math. J. {\bf 84}, 217--235, (1996).}\\

\mbox{[S]} & \mbox{Sorger C., {\it La formule de Verlinde}, S\'eminare Bourbaki {1994/1995.}}\\

\mbox{[TUY]} & \mbox{Tsuchiya A., Ueno K. and Yamada Y., {\it Conformal field theory on universal}}\\
            &\mbox{{\it family of stable curves with gauge symmetries}, Adv. Stud. Pure Math. {\bf 19},}\\
            &\mbox{459--565, (1989).}\\

\mbox{[TW]} &\mbox{Teleman C. and Woodward C., {\it Parabolic bundles, products of conjugacy}}\\
    &\mbox{{\it classes, and Gromov-Witten invariants}, Annales de l'institut Fourier, {\bf 53}, 3,}\\
    &\mbox{(2003), 713--748}\\

\mbox{[V]} & \mbox{Verlinde E., {\it Fusion rules and modular transformations in 2d conformal field }}\\
            &\mbox{{\it theory,} Nucl. Phys. B {\bf 300}, 360--376, (1988).}

\end{array} \]

\vskip 2pt
\noindent International Centre for Theoretical Physics (ICTP), Mathematics Section\\
\noindent Strada Costiera 11, I-34014 Trieste, Italy.\\
\noindent {\tt aboysal}@{\tt ictp.it}\\
\end{document}